\documentclass[12pt,a4paper]{article}
\usepackage{amsmath}
\usepackage{amssymb}
\usepackage{amsthm}
\usepackage{amsfonts}
\usepackage{latexsym}

\theoremstyle{plain}
\newtheorem{thm}{Theorem}[section]
\newtheorem{cor}{Corollary}[section]
\newtheorem{lem}{Lemma}[section]

\theoremstyle{definition}
\newtheorem{defn}{Definition}[section]

\theoremstyle{remark}

\newtheorem{rem}{Remark}[section]

\title{On the Weyl law for Toeplitz operators}
\author{Roberto Paoletti\footnote{\noindent{\bf Address:}
Dipartimento di Matematica e Applicazioni, Universit\`a degli Studi
di Milano Bicocca, Via R. Cozzi 53, 20125 Milano,
Italy; {\bf e-mail}: roberto.paoletti@unimib.it }}
\date{}

\begin{document}

\maketitle

\section{Introduction}

A  Weyl law for Toeplitz operators was proved
by Boutet de Monvel and Guillemin  (\cite{bg}, Theorem
13.1), using the theory of Fourier-Hermite distributions and their symbolic calculus based on
symplectic spinors. The result holds in the setting of so-called Toeplitz structures, a generalization
of strictly pseudoconvex domains, and the operators involved act on certain spaces
of half-forms. The circle bundle of a positive line bundle on a symplectic manifold, with an appropriate
connection, is an important
special case of a Toeplitz structure.

On the other hand, a circle of ideas has emerged recently
(\cite{z}, \cite{bsz}, \cite{sz}), in which asymptotic properties of positive line bundles
are studied by appealing directly to the microlocal description of Szeg\"{o} kernels
in \cite{bs}, and to its generalization to the almost complex symplectic setting in \cite{sz}.
These ideas and techniques lead to
several rather explicit local and global
asymptotic expansions related to the (generalized) Szeg\"{o} kernel of the line
bundle. It seems natural to ask whether the Weyl law can also be interpreted in light of this approach.

Accordingly, the aim of this paper is to provide an alternative,
perhaps relatively elementary proof of the
theorem of Boutet de Monvel and Guillemin, in the case of Toeplitz structures related to positive line bundles
on symplectic manifolds. As is well known, this setting is natural in algebraic and symplectic geometry,
and in the theory of geometric quantization.
For ease of exposition, we shall state our results in the context of
complex projective geometry, but in view of the theory of
\cite{sz} the present arguments apply to the more general
almost complex symplectic category.

In short, by pairing techniques from \cite{sz} with certain classical arguments used to estimate the spectral function
of a pseudodifferential operator \cite{h}, \cite{gs}, we shall obtain local Weyl
laws for the spectral function of a Toeplitz operator; the global Weyl law then follows by integration.
Unlike \cite{bg}, the operators will be understood to act on half-densities, rather that half-forms.

Before stating the result in point,
it is in order to fix some conventions, and to recall some basic notions from
\cite{bg}.

Suppose $M$ is a d-dimensional complex projective manifold, and let $A$ be an ample line
bundle on $M$. There exists an Hermitian metric on $A$, such that the
unique compatible connection has normalized curvature $\Theta=-2i\,\omega$, where
$\omega$ is a K\"{a}hler form. The induced volume form $dV_M=\frac{1}{\mathrm{d}!}\,\omega^{\wedge \mathrm{d}}$
on $M$ has
total volume $\mathrm{vol}(M)=\frac{\pi^\mathrm{d}}{\mathrm{d}!}\,\int _Mc_1(A)^\mathrm{d}$.

With this assumption, the unit disc bundle $D\subseteq A^*$ is a strictly
pseudoconvex domain, with boundary $X=\partial D$ the unit circle bundle in $A^*$.
The compatible connection on $A$ determines
the connection form $\alpha\in \Omega^1(X)$, and by positivity
$(X,\alpha)$ is a contact manifold. In particular, $X$ has the volume form
$\frac{1}{2\pi}\,dV_X=\frac{1}{2\pi}\,\alpha\wedge \pi^*(dV_M)$;
given this, we shall identify densities, half-densities, and functions on $X$.

Given the $L^2$-product on $X$ induced by
$d\mu_X=:\frac{1}{2\pi}\,dV_X$,
let $H(X)\subseteq L^2(X)$ denote  the
Hardy space of $X$, and $\Pi:L^2(X)\rightarrow H(X)$ the orthogonal projector. The Szeg\"{o} kernel of $X$ is the
distributional kernel of $\Pi$, which we shall denote by the same symbol, $\Pi\in \mathcal{D}'(X\times X)$.
If $\{s_j\}$ is an orthonormal basis of $H(X)$, then $\Pi(x,y)=\sum_{j=1}^{+\infty}s_j(x)\,\overline{s_j(y)}$ and
\begin{equation}
\label{eqn:pairing-density}
\Pi(f)(x)=\int _X\Pi(x,y)\,f(y)d\mu_X(y)\,\,\,\,\,\,\,\,\,\,\left(f\in L^2(X),\,x\in X\right).
\end{equation}

We shall deal with the following class of operators:

\begin{defn}
\label{defn:toeplitz}
If $m\in \mathbb{Z}$, a Toeplitz operator of order $m$ on $X$ is an operator
$T:\mathcal{D}'(X)\rightarrow \mathcal{D}'(X)$
of the form
$T=:\Pi\circ P\circ \Pi$, where $P$ is a pseudodifferential operator of classical type
and order $m$ on $X$.
\end{defn}

The composition in Definition \ref{defn:toeplitz}
is well-defined,
in view of the microlocal structure of $\Pi$ \cite{bs} and standard results on wave fronts \cite{d}.

Now $\alpha$ being a connection form on $X$ is equivalent to the condition that
the closed cone
\begin{equation}
\label{eqn:sigma-alpha}
\Sigma=:\left\{(x,r\alpha_x):x\in X,\,r>0\right\}\subseteq T^*X\setminus \{0\}
\end{equation}
be a symplectic submanifold of the cotangent bundle of $X$, equipped with its canonical symplectic structure.
By the theory in \S 1 and \S 2 of \cite{bg},
the following definition is well-posed:

\begin{defn}
\label{defn:toeplitz-symb}
Let $T:\mathcal{D}'(X)\rightarrow \mathcal{D}'(X)$ be a Toeplitz operator of order $m$,
and let $\Sigma$ be as in (\ref{eqn:sigma-alpha}).
If $P$ is as in Definition
\ref{defn:toeplitz}, the \textit{symbol} $\sigma_T$ of $T$ is the $\mathcal{C}^\infty$ function
on $\Sigma$ given by the restriction of the symbol of $P$.
Furthermore, $T$ is \textit{elliptic} if $\sigma_T$ is the restriction to $\Sigma$
of an elliptic symbol of degree $m$.
\end{defn}

Thus a Toeplitz operator $T$ determines a $\mathcal{C}^\infty$-function
$\varsigma_T$ on $X$, defined by
$\varsigma_T(x)=:\sigma_T(x,\alpha_x)$ ($x\in X$).

We shall say that the Toeplitz operator $T$ is self-adjoint
to mean that it is formally self-adjoint with respect to
the $L^2$-product associated to $d\mu_X$; in this case, $T$ is
an essentially self-adjoint operator with domain $\mathcal{C}^{\infty}(X)$,
and its symbol $\sigma_T$ is real valued.

Suppose that $T$ is a self-adjoint and elliptic Toeplitz operator
of degree $1$, with
$\sigma_T>0$. By the theory in \S 1 and \S 2 of \cite{bg},
$T\ge -C\,I$ for some constant $C\in \mathbb{R}$.
Furthermore, its spectrum $\mathrm{Spec}(T)$ is discrete, bounded from below and has
only $+\infty$ as an accumulation point (Proposition 2.14 of \cite{bg}); if we view $T$ as
an operator $T_H$ with domain $\mathcal{C}^{\infty}(X)\cap H(X)$,
every eigenvalue occurs with finite multiplicity.

Thus $\mathrm{Spec}\left(T_H \right)=\big\{\lambda_j\big\}\subseteq \mathbb{R}$, with $\lambda_1\le \lambda_2\le \cdots$,
where each $\lambda_j$ is repeated according to its multiplicity, and $\lambda_j\rightarrow +\infty$ for $j\rightarrow
+\infty$. For every $j=1,2,\ldots$, let $e_j\in H (X)\cap \mathcal{C}^\infty(X)$ be an eigenfunction of
$T$ corresponding to the eigenvalue $\lambda_j$ and of unit $L^2$-norm.
We can arrange that $\{e_j\}$ is an orthonormal basis of $H(X)$.

For every $\lambda\in \mathbb{R}$,
let us set
$$
\mathbb{T}_\lambda=:\mathrm{span}\Big\{e_j:\lambda_j\le \lambda\Big\},
$$
a finite dimensional subspace of $H(X)$, and denote by
$
\mathfrak{T}_\lambda :L^2(X)\rightarrow \mathbb{T}_\lambda
$ the orthogonal projector:
$\mathfrak{T}_\lambda$ is the smoothing operator with distributional kernel
the \textit{spectral function}
$$
\mathcal{T}(\lambda,x,y)=:\sum _{\lambda_j\le \lambda}e_j(x)\cdot \overline{e_j(y)}\,\,\,\,\,\,\,\,\,\,\,\,
(x,y\in X).
$$
We shall also let $\mathcal{T}(\lambda,x)=:\mathcal{T}(\lambda,x,x)$ ($x\in X,\,\lambda\in \mathbb{R}$).

\begin{thm}
\label{thm:local-weyl-law}
Let $T$ be a self-adjoint first order elliptic Toeplitz operator on $X$ with $\sigma_T>0$.
Then as $\lambda\rightarrow +\infty$ we have
$$
\mathcal{T}(\lambda,x)=\frac{\pi}{d+1}\cdot\left(\frac{\lambda}{\pi\, \varsigma_T(x)}\right)
^{\mathrm{d}+1}+O\left(\lambda^\mathrm{d}\right),
$$
uniformly in $x\in X$.
\end{thm}

The \textit{counting function} of $T$ is
\begin{equation}
\label{eqn:counting-fctn}
N_T(\lambda)=:\sharp \Big\{j\in \mathbb{N}\,:\,\lambda_j\le \lambda\Big\}\,\,\,\,\,\,\,\,\,\,\,\,
(\lambda\in \mathbb{R}).
\end{equation}
Clearly,
\begin{equation}
\label{eqn:counting-trace}
N_T(\lambda)=\mathrm{trace}(\mathfrak{T}_\lambda)=\int _X \mathcal{T}(\lambda,x)\,\mathrm{vol}_X(x).
\end{equation}

To check our normalizations, suppose
$T=\Pi\circ \left(\frac{1}{i}\,\frac{\partial}{\partial \theta}\right)\circ \Pi$, where $\frac{\partial}{\partial \theta}$
denotes the generator of the $S^1$-action.
Then by the Riemann-Roch Theorem
\begin{equation}
\label{eqn:riemann-roch}
N_T(\lambda)=\sum _{k\le \lambda}h^0\left (M,A^{\otimes k}\right)=\frac{1}{(\mathrm{d}+1)!}\,
\lambda^{\mathrm{d}+1}\,\int _M c_1(A)^\mathrm{d}+O\left(\lambda^\mathrm{d}\right)
\end{equation}
for $\lambda\rightarrow +\infty$.
Now Theorem \ref{thm:local-weyl-law} says that
$\mathcal{T}(\lambda,x)=\frac{\pi}{d+1}\cdot\left(\frac{\lambda}{\pi}\right)
^{\mathrm{d}+1}+O\left(\lambda^\mathrm{d}\right)$, and therefore predicts the estimate
$$N_T(\lambda)=\frac{\pi}{d+1}\cdot\left(\frac{\lambda}{\pi}\right)
^{\mathrm{d}+1}\cdot \mathrm{vol}(M)+O\left(\lambda^\mathrm{d}\right),
$$
in agreement with (\ref{eqn:riemann-roch}).

More generally, the Weyl law of \cite{bg} follows by using Theorem
\ref{thm:local-weyl-law} in (\ref{eqn:counting-trace}):
\begin{cor}
\label{cor:global-weyl-law}
Let
$\Sigma_1\subseteq \Sigma$ be the locus where $\sigma_T<1$.
In the situation of Theorem \ref{thm:local-weyl-law},
as $\lambda\rightarrow +\infty$ we have
$$
N_T(\lambda)=
\left(\frac{\lambda}{2\pi}\right)^{\mathrm{d}+1}\,\mathrm{vol}(\Sigma_1)+O\left(\lambda^\mathrm{d}\right).
$$
\end{cor}

\section{Proof of Theorem \ref{thm:local-weyl-law}.}

By Lemma 12.1 of \cite{bg}, there exists a
self-adjoint first-order pseudodifferential operator of classical type
$Q$ on $X$ such that
$[Q,\Pi]=0$, $q>0$ everywhere, and $T=\Pi\circ Q\circ \Pi$; here
$q:T^*X\setminus \{0\}\rightarrow \mathbb{R}$ is the symbol of
$Q$.

Let us consider the one-parameter group of $L^2$-unitary operators
$$
U(\tau)=:e^{i\tau Q}\,\,\,\,\,\,\,(\tau \in \mathbb{R})
$$
(\S 12 of \cite{gs}, \S 12 of \cite{bg}).
%Furthermore, $U(\tau)$ has the following microlocal description.

The $\mathcal{C}^\infty$-function $q$
generates a one-parameter group $\Phi_\tau$ ($\tau\in \mathbb{R}$)
of homogeneous symplectomorphisms
of $T^*X\setminus \{0\}$. For every $\tau\in \mathbb{R}$
the graph of $\Phi_\tau$, $\mathrm{graph}(\Phi_\tau)$, is therefore
a conic Lagrangian relation, and
$U(\tau)$ is a Fourier integral operator of degree $0$ associated to it.
In fact, there exists $\epsilon_0>0$ such that if $|\tau|<\epsilon_0$
then $U(\tau)=V(\tau)+R(\tau)$, where $V(\tau)$ and $R(\tau)$ are as follows.

First, $R(\tau)$ is a smoothing operator, and in fact
its kernel $R(\tau,x,y)$ is a $\mathcal{C}^\infty$ function on
$\big((-\epsilon_0,\epsilon_0)\times X\times X\big)$ (Proposition 12.3 of \cite{gs}).

Next,
$V(\tau)$ can be locally
represented by a Fourier integral in the following manner.
On a coordinate patch, with abuse of notation let us identify $x,y,\cdots\in X$
with their local coordinates in $\mathbb{R}^{2\mathrm{d}+1}$.
Let
$\varphi (\tau,\cdot,\cdot)$ be a generating function
for $\Phi_\tau$. Since $\Phi_0$ is the identity and $\varphi$ satisfies the Hamilton-Jacobi equation,
\begin{equation}
\label{eqn:generating-for-small-tau}
\varphi(\tau,x,\eta)=x\cdot \eta +\tau \,q(x,\eta)+O\left(\tau^2\right)\cdot \|\eta\|.
\end{equation}
Then
\begin{equation}
\label{eqn:local-representation-V}
V(\tau)(x,y)=\frac{1}{(2\pi)^{2\mathrm{d}+1}}\,
\int _{\mathbb{R}^{2\mathrm{d}+1}}e^{i[\varphi(\tau,x,\eta)-y\cdot \eta]}\,a(\tau,x,y,\eta)\,d\eta.
\end{equation}
Here $a(\tau,\cdot,\cdot,\cdot)\in S^0_{\mathrm{cl}}$
for every $\tau$. To determine $a(0,\cdot,\cdot,\cdot)$, remark that
$U(0)=\mathrm{Id}$; therefore, if $V(\tau)$ was to act on functions by means of the Euclidean density $dx$
on $\mathbb{R}^{2\mathrm{d}+1}$, we would have $a(0,\cdot,\cdot,\cdot)=1$. Since however $V(\tau)$ acts
on functions by means of the density $d\mu_X$ (that is, as in (\ref{eqn:pairing-density}) with $V(\tau)$
in place of $\Pi$), we must have $a(0,\cdot,\cdot,\cdot)=1/\mathcal{V}(y)$, where
$\mathcal{V}\,dx$ is the local coordinate expression of $d\mu_X$.

By \cite{bs}, $\Pi$ is also a Fourier integral operator, with complex phase however. More precisely,
 $\Pi=\Pi _1+\Pi_2$, where $\Pi_2$ is smoothing, while locally
\begin{equation}
\label{eqn:local-representation-Pi}
\Pi_1(x,y)=\int _0^{+\infty}e^{it\psi(x,y)}\,s(t,x,y)\,dt,
\end{equation}
where $\Im \psi\ge 0$ and $s(t,x,y)\sim \sum _{j=0}^{+\infty}t^{\mathrm{d}-j}\,s_j(x,y)$.

If $\epsilon>0$ and
$\chi\in \mathcal{C}^\infty _0\big((-\epsilon,\epsilon)\big)$, for any suitable family of operators
$A(\tau)$ we set (see \S 12 of \cite{gs})
\begin{equation*}
\label{eqn:A-chi}
A_\chi=:\int _{-\epsilon}^\epsilon \chi(\tau)\,A(\tau)\,d\tau.
\end{equation*}

Let us pick a sufficiently small $\epsilon \in (0,\epsilon_0)$, and apply this to
$S(\tau)=:\Pi\circ U(\tau)\circ \Pi$ ($\tau\in \mathbb{R}$).
By the hypothesis, actually $S(\tau)=U(\tau)\circ \Pi$.
Thus,
\begin{equation}
\label{eqn:s-chi}
S_\chi=U_\chi\circ \Pi=V_\chi\circ \Pi+R_\chi\circ \Pi.
\end{equation}
Since $U_\chi$ is a smoothing operator \cite{gs}, so is $S_\chi$.

When $\chi$ is replaced by $\chi \cdot e^{-i\lambda \,(\cdot)}$, and (\ref{eqn:s-chi})
is rewritten with some abuse of language in terms of distributional kernels, we obtain
\begin{eqnarray}
\label{eqn:s-chi-kernel}
S_{\chi \,e^{-i\lambda \,(\cdot)}} (x,y)&=&\left( U_{\chi \,e^{-i\lambda \,(\cdot)}}\circ \Pi\right)(x,y)\\
&=&\left( V_{\chi \,e^{-i\lambda \,(\cdot)}}\circ \Pi\right)(x,y)
+\left( R_{\chi \,e^{-i\lambda \,(\cdot)}}\circ \Pi\right)(x,y)
\,\,\,\,\,\,\,(x,y\in X).\nonumber
\end{eqnarray}
We shall now study the asymptotics of (\ref{eqn:s-chi-kernel}) for $\lambda\rightarrow \infty$.

To begin with, the kernel $\mathcal{R}(t,x,y)$ of $R(t)\circ \Pi$ is a $\mathcal{C}^\infty$ function on
$(-\epsilon_0,\epsilon_0)\times X\times X$.
The second summand in (\ref{eqn:s-chi}) is then given by
\begin{equation}
\label{eqn:rapid-decay}
\int _{-\epsilon}^\epsilon e^{-i\lambda t}\,\chi(t)\,\mathcal{R}(t,x,y)\,dt=\widehat{\chi\,\mathcal{R}}(\lambda,x,y)=
O\left (\lambda ^{-\infty}\right)
\end{equation}
for $\lambda\rightarrow \infty$, uniformly in $x,y\in X$;
here $\,\widehat{}\,$ is the
Fourier transform with respect to $t$.

Let $\sim$ denote equality of asymptotic expansions.
By (\ref{eqn:s-chi}) and (\ref{eqn:rapid-decay}), for all $x\in X$
\begin{eqnarray}
\label{eqn:s-chi-lambda-diag}
\lefteqn{
S_{\chi \,e^{-i\lambda \,(\cdot)}} (x,x)\sim \left(V_{\chi \,e^{-i\lambda \,(\cdot)}}\circ \Pi\right)(x,x)}\\&=&\int _X\int_{-\epsilon}^\epsilon e^{-i\lambda \tau}\,\chi (\tau)\,
V(\tau)(x,z)\,\Pi(z,x)\,d\mu_X(z)\,d\tau.                                                                           \nonumber
\end{eqnarray}

Let us fix $x\in X$ and an arbitrarily small open neighborhood
$X_1\subseteq X$; let furthermore
$X_2\subseteq X$ be an open neighborhood of $X\setminus X_1$, such that $x\not\in \overline{X_2}$.
Let $\big\{\varrho_j\big\}$ be a a partition of unity
on $X$ subordinate to the open cover $\mathcal{U}=:\left\{X_1,X_2\right\}$ of $X$.
We obtain from (\ref{eqn:s-chi-lambda-diag}):
\begin{eqnarray}
\label{eqn:s-chi-lambda-diag-part}
\lefteqn{
S_{\chi \,e^{-i\lambda \,(\cdot)}} (x,x)\sim \left( V_{\chi \,e^{-i\lambda \,(\cdot)}}\circ \Pi\right)(x,x)}\\
&=&\sum _{j=1}^2\int _{X_j}\int_{-\epsilon}^\epsilon e^{-i\lambda \tau}\,\chi (\tau)\,\varrho _j(z)\,
V(\tau)(x,z)\Pi(z,x)\,d\mu_X(z)\,d\tau.                                                                           \nonumber
\end{eqnarray}

\begin{lem}
\label{lem:rapid-decay-second-term}
As $\lambda\rightarrow \infty$, we have:
$$
\int _{X_2}\int_{-\epsilon}^\epsilon e^{-i\lambda \tau}\,\chi (\tau)\,\varrho _2(z)
V(\tau)(x,z)\,\Pi(z,x)\,d\mu_X(z)\,d\tau=O\left(\lambda^{-\infty}\right).
$$
\end{lem}

\textit{Proof.} Since the singular support of
$\Pi$ is the diagonal in $X\times X$,
$\Pi(\cdot,x)$ is
$\mathcal{C}^\infty$ on $X_2$.
If $p_x(z)=:\varrho _2(z)\,\Pi(z,x)$, we have $p_x\in \mathcal{C}^\infty_0(X_2)$.

Now $V(\tau)$ is a $\mathcal{C}^\infty$ family of distributions on $X\times X$,
and by regularity
$V(\tau,x)=:V(\tau)(x,\cdot)$ is therefore a  well-defined $\mathcal{C}^\infty$-family of distributions on $X$,
parametrized by $\mathbb{R}\times X$.

Consequently,
$\gamma (\tau)=:\chi(\tau)\,\langle V(\tau,x),p_x\rangle$ is a compactly supported
$\mathcal{C}^\infty$ function of $\tau$, hence
its Fourier transform is rapidly decreasing.

\hfill Q.E.D.

\medskip

\begin{rem}
If $X_2$ is defined by an inequality on the distance
function from $x$ (using the natural metric on $X$), and the same
bound is used for every $x\in X$, the previous statement
holds uniformly in $x$.
\end{rem}

Recalling that $\Pi=\Pi_1+\Pi_2$, we conclude that

\begin{eqnarray}
\label{eqn:s-chi-lambda-diag-part-reduced}
\lefteqn{S_{\chi \,e^{-i\lambda \,(\cdot)}} (x,x)\sim }\\&=&\sum _{j=1}^2\int _{X_1}
\int_{-\epsilon}^\epsilon e^{-i\lambda \tau}\,\chi (\tau)\,\varrho _1(z)
V(\tau)(x,z)\,\Pi_j(z,x)\,d\mu_X(z)\,d\tau.\nonumber
\end{eqnarray}

Since $\Pi_2$ is smoothing, a similar argument to the proof of Lemma \ref{lem:rapid-decay-second-term}
shows that the second summand in (\ref{eqn:s-chi-lambda-diag-part-reduced}) is $O\left(\lambda^{-\infty}\right)$
as $\lambda\rightarrow \infty$.
To estimate asymptotically the first summand, we make use of (\ref{eqn:generating-for-small-tau}),
(\ref{eqn:local-representation-V}) and (\ref{eqn:local-representation-Pi}) and obtain

\begin{eqnarray}
\label{eqn:s-chi-lambda-diag-part-reduced-P-2-solo-fourier}
\lefteqn{S_{\chi \,e^{-i\lambda \,(\cdot)}} (x,x)\sim }\\
&=&\frac{1}{(2\pi)^{2\mathrm{d}+1}}\,\int _{X_1}\int _0^{+\infty}\int_{-\epsilon}^\epsilon \int_{\mathbb{R}^{2\mathrm{d}+1}}
e^{i\,\Phi_1(x,z,t,\tau,\eta,\lambda)}\,A(x,z,t,\tau,\eta)\,\,d\mu_X(z)\,dt\,d\tau\,d\eta, \nonumber
\end{eqnarray}
where
\begin{eqnarray}
\label{eqn:phase-A}
\Phi_1(x,z,t,\tau,\eta,\lambda)&=:&\varphi (\tau,x,\eta)-z\cdot \eta+t\psi(z,x)-\lambda \tau\\
&=&
(x-z)\cdot \eta+\tau\,q(x,\eta)+t\psi(z,x) -\lambda \tau+O\left(\tau^2\right)\cdot \|\eta\|,\nonumber
\end{eqnarray}
and
\begin{eqnarray}
\label{eqn:amplitude-A}
A(x,z,t,\tau,\eta)=:\chi (\tau)\,\varrho_1(z)\,a(\tau,x,z,\eta)\,s(t,z,x).
\end{eqnarray}

Since $q$ is positively homogenous of degree one in $\eta$ and everywhere positive for
$\eta\neq \mathbf{0}$, (\ref{eqn:phase-A}) implies that
$\partial _\tau\Phi\sim \|\eta\|+|\lambda|$ for $\lambda\rightarrow -\infty$; integrating by parts in $\tau$,
$S_{\chi \,e^{-i\lambda \,(\cdot)}} (x,x)=O\left(\lambda^{-\infty}\right)$ for $\lambda\rightarrow
-\infty$.

To study the asymptotics of (\ref{eqn:s-chi-lambda-diag-part-reduced-P-2-solo-fourier})
for $\lambda\rightarrow +\infty$, following \cite{gs} we choose a radial function
$F\in \mathcal{C}^\infty _0\left(\mathbb{R}^{2\mathrm{d}+1}\setminus \{\mathbf{0}\}\right)$
identically equal to $1$ for $\frac 1C\le \|\eta \|\le C$ for some $C\gg 0$,
and split (\ref{eqn:s-chi-lambda-diag-part-reduced-P-2-solo-fourier}) as:
\begin{eqnarray}
\label{eqn:s-chi-lambda-diag-part-reduced-P-2-solo-fourier-spezzato}
\lefteqn{S_{\chi \,e^{-i\lambda \,(\cdot)}} (x,x)\sim }\\
&=&\frac{1}{(2\pi)^{2\mathrm{d}+1}}\,\left\{\int _{X_1}\int _0^{+\infty}\int_{-\epsilon}^\epsilon \int_{\mathbb{R}^{2\mathrm{d}+1}}
e^{i\,\Phi_1(x,z,t,\tau,\eta,\lambda)}\,A(x,z,t,\tau,\eta)\,F\left (\frac \eta \lambda \right)\,d\mu_X(z)\,dt\,d\tau\,d\eta \right.\nonumber \\
&&+\left.\int _{X_1}\int _0^{+\infty}\int_{-\epsilon}^\epsilon \int_{\mathbb{R}^{2\mathrm{d}+1}}
e^{i\,\Phi_1(x,z,t,\tau,\eta,\lambda)}\,A(x,z,t,\tau,\eta)\,\left[1-F\left (\frac \eta \lambda \right)\right]\,d\mu_X(z)\,dt\,d\tau\,d\eta\right\}.\nonumber
\end{eqnarray}
Since $\|\partial _\tau \Phi\|\sim
\|\eta\|+\lambda$ where $1-F\left (\frac \eta \lambda \right)\neq 0$, the second term in
(\ref{eqn:s-chi-lambda-diag-part-reduced-P-2-solo-fourier-spezzato}) is
$O\left(\lambda ^{-\infty}\right)$ for $\lambda\rightarrow +\infty$.

To estimate the first term in (\ref{eqn:s-chi-lambda-diag-part-reduced-P-2-solo-fourier-spezzato}),
let us
perform the change of variables $t\rightsquigarrow \lambda t$,
$\eta \rightsquigarrow \lambda \eta$, and use the homogeneity of $\Phi$ to obtain:
\begin{eqnarray}
\label{eqn:s-chi-lambda-diag-part-reduced-P-2-solo-fourier-riscalato}
S_{\chi \,e^{-i\lambda \,(\cdot)}} (x,x)
&\sim&\frac{\lambda^{2\mathrm{d}+2}}{(2\pi)^{2\mathrm{d}+1}}\,\int _{X_1}\int _0^{+\infty}\int_{-\epsilon}^\epsilon \int _{1/D}^D
\int_{S^{2\mathrm{d}}}
e^{i\,\lambda\,\Phi_1(x,z,t,\tau,r\omega,1)}\,A(x,z,\lambda t,\tau,\lambda r\omega)\nonumber\\
&&\cdot r^{2\mathrm{d}}\,F(r)\,
d\mu_X(z)\,dt\,d\tau\,dr\,d\omega,
\end{eqnarray}
where now integration in $dr$ takes place over the interval
$\frac 1D\le r\le D$ for some $D\gg 0$.

To proceed further, we shall need a more specific choice of
a system of local coordinates near $x$. It
will be useful to work in a system of \textit{Heisenberg local coordinates}
centered at $x$, as defined in \cite{sz}, for in these
systems of coordinates scaling limits of Szeg\"{o} kernel
exhibit their universal nature (the focus in \cite{sz} is
on the asymptotics for $k\rightarrow +\infty$, where $k$ indexes the
isotype with respect to the $S^1$-action, but we shall exploit some
formal analogies with the present computations). We refer to
\cite{sz} for a definition of Heisenberg local
coordinates.

Suppose then that $(\theta,\mathrm{v}):X_1\rightarrow (-a,a)\times B_{2\mathrm{d}}\big(\mathbf{0},\delta\big)$
is a system of Heisenber local coordinates centered at $x$, for some $a,\delta>0$; here $B_{2\mathrm{d}}\big(\mathbf{0},\delta\big)\subseteq \mathbb{R}^{2\mathrm{d}}
\cong \mathbb{C}^\mathbf{d}$ is the open ball of radius $\delta$ centered at the origin.
Following \cite{sz},
let $z=x+(\theta,\mathbf{v})$ denote the point with coordinates $(\theta,\mathbf{v})$.
Furthermore, we shall write
$\omega =(\omega_0,\omega_1)$, where $\omega_0\in \mathbb{R}$, $\omega_1\in \mathbb{R}^{2\mathrm{d}}$ and
$\omega_0^2+\|\omega_1\|^2=1$.
On the upshot, the phase in (\ref{eqn:s-chi-lambda-diag-part-reduced-P-2-solo-fourier-riscalato})
becomes
\begin{eqnarray}
\label{eqn:phase-phi-heis}
\lefteqn{\Phi_2\big(x,\theta,\mathbf{v},t,\tau,r,\omega\big)=:\Phi_1(x,z,t,\tau,r\omega,1)}\\
&=:&-r\,\mathbf{v}\cdot\omega_1
+t\psi\Big(x+\left(\theta,\mathbf{v}\right),x\Big)-r\,\theta \,\omega_0
+r\tau\,q\big(x,\omega\big)-\tau+O\left(\tau^2\right)\,r.
\nonumber
\end{eqnarray}
The paring $\mathbf{v}\cdot\omega_1$ is the standard scalar product as vectors in $\mathbb{R}^{2\mathrm{d}}$.

Let $d\mu_X(z)=\mathcal{V}(\theta,\mathbf{v})\,d\theta\,d\mathbf{v}$
be the volume density expressed in local coordinates.
Integration in $d\mathbf{v}$ takes place over the open ball
$B_{2\mathrm{d}}(\mathbf{0},\delta)\subseteq \mathbb{C}^\mathrm{d}$.
By definition of Heisenberg local coordinates,
$\mathcal{V}(\theta,\mathbf{0})=1/2\pi$ for every $\theta$.

We next show that, up to a rapidly decreasing contribution, the integration in $d\omega$ can be
restricted to an arbitrarily small open neighborhood of $(1,\mathbf{0})$ in $S^{2\mathrm{d}}$.
Suppose $0<\delta_1<\delta_2<1$, and define
$$S_-=:\left\{\omega \in S^{2\mathrm{d}}:\omega _0<\delta_2\right\},\,\,\,\,\,
S_+=:\left\{\omega \in S^{2\mathrm{d}}:\omega _0>\delta_1\right\}.
$$
Let $\{\beta_-,\beta_+\}$ be a partition of unity subordinate to
the open cover $\mathcal{U}=\left\{S_-,S_+\right\}$
of $S^{2\mathrm{d}}$.
Inserting the identity $\beta_-+\beta_+=1$ in
(\ref{eqn:s-chi-lambda-diag-part-reduced-P-2-solo-fourier-riscalato}), we obtain
$$
S_{\chi \,e^{-i\lambda \,(\cdot)}} (x,x)\sim S_{\chi \,e^{-i\lambda \,(\cdot)}} (x,x)_-+
S_{\chi \,e^{-i\lambda \,(\cdot)}} (x,x)_+,
$$
where in the former term the integrand in (\ref{eqn:s-chi-lambda-diag-part-reduced-P-2-solo-fourier-riscalato})
has been multiplied by $\beta_-(\omega)$, in the latter by $\beta_+(\omega)$.

\begin{lem}
\label{lem:reduction-on-sphere}
If $\epsilon$ is sufficiently small,
then $S_{\chi \,e^{-i\lambda \,(\cdot)}} (x,x)_-=O\left(\lambda ^{-\infty}\right)$ when
$\lambda\rightarrow +\infty$.
\end{lem}

\textit{Proof.}
If in the integral
representation (\ref{eqn:local-representation-V}) for $V(\tau)$ a cut-off is introduced in the neighborhood of
$\eta=0$, only a smoothing term is lost. The argument leading to (\ref{eqn:rapid-decay})
implies that the corresponding contribution to $S_{\chi \,e^{-i\lambda \,(\cdot)}} (x,x)$
is $O\left(\lambda^{-\infty}\right)$.

With this cut-off implicit, the identity
$\beta_-\big(\eta/\|\eta\|\big)+\beta_+\big(\eta/\|\eta\|\big)=1$ in (\ref{eqn:local-representation-V})
leads to
a decomposition
$V(\tau)=V(\tau)_-+V(\tau)_+$.
Clearly, $S_{\chi \,e^{-i\lambda \,(\cdot)}} (x,x)_-$ is
obtained by replacing $V(\tau)$ by $V(\tau)_-$ in the previous construction.

The wave front of $V(\tau)$ is locally parametrized by
$$
\left(x',\eta\right)\mapsto \left(x',\frac{\partial \varphi}{\partial x}\left(\tau,x',\eta\right),
\frac{\partial \varphi}{\partial \eta}\left(\tau,x',\eta\right),-\eta\right).
$$
Let us write $\eta=\big(\eta_0,\eta_1\big)\in \mathbb{R}^{2\mathrm{d}+1}
\cong\mathbb{R}\times \mathbb{R}^{2\mathrm{d}}$.
By construction the wave front of $V(\tau)_-$ satisfies
\begin{eqnarray}
\label{eqn:wave-front-of-V-a}
\lefteqn{\mathrm{WF}\left(V(\tau)_-\right)}\\
&\subseteq &\left\{\left(x',\frac{\partial \varphi}{\partial x}\left(\tau,x',\eta\right),
\frac{\partial \varphi}{\partial \eta}\left(\tau,x',\eta\right),-\eta\right):x'\in X_1,\,
\eta\in \mathbb{R}^{2\mathrm{d}+1}:
\eta_0/\|\eta\|<\delta_2\right\}\nonumber\\
&=&\left\{\left(\Phi_\tau\left(\frac{\partial \varphi}{\partial \eta}\left(\tau,x',\eta\right),\eta\right),
\left(\frac{\partial \varphi}{\partial \eta}\left(\tau,x',\eta\right),-\eta\right)\right):x'\in X_1,\,
\eta\in \mathbb{R}^{2\mathrm{d}+1}:
\eta_0/\|\eta\|<\delta_2\right\}.\nonumber
\end{eqnarray}
On the other hand, the wave front of $\Pi$ is
\begin{equation}
\label{eqn:wave-front-of-Pi}
  \mathrm{WF}(\Pi) =\big\{(x,r\alpha_x,x,-r\alpha_x):x\in X,\,r>0\big\}.
\end{equation}

By definition of Heisenberg local coordinates, the cotangent vector $(x,\alpha_x)$
corresponds to $\Big((0,\mathbf{0}),(1,\mathbf{0})\Big)$.
Therefore, if $X_1$ has been chosen sufficiently small, there exists $a>0$
such that (mixing intrinsic notation and expressions in local coordinates)
$\mathrm{dist}_X\left(\left (\frac{\partial \varphi}{\partial \eta}\left(\tau,x',\eta\right),\eta\right),(x,\alpha_x)\right)>a$,
if $x\in X_1$, $\|\eta\|=1$ and
$\eta_0<\delta_2$, where $\mathrm{dist}_X$ is the Riemannian distance on $X$.
Therefore, if $\epsilon>0$ is sufficiently small we also have
$\mathrm{dist}_X\left(\Phi_\tau\left (\frac{\partial \varphi}{\partial \eta}\left(\tau,x',\eta\right),\eta\right),(x,\alpha_x)\right)>\frac a2$ in the same range if $|\tau|<\epsilon$.

Given this, (\ref{eqn:wave-front-of-V-a}), and (\ref{eqn:wave-front-of-Pi}) we have
$\mathrm{WF}'\big(V(\tau)_-\circ \Pi\big)=\mathrm{WF}'\big(V(\tau)_-\big)\circ \mathrm{WF}'\big(\Pi\big)=\emptyset$;
here $\mathrm{WF}'$ is the image of $\mathrm{WF}$ under the map
$(x,\xi,y,\eta)\mapsto (x,\xi,y,-\eta)$. Hence, $V(\tau)_-\circ \Pi$
is a smoothing operator for every $\tau\in (-\epsilon,\epsilon)$. The claim
follows by arguing again as for (\ref{eqn:rapid-decay}).

\hfill Q.E.D.

\medskip

Given (\ref{eqn:s-chi-lambda-diag-part-reduced-P-2-solo-fourier-riscalato})
and Lemma \ref{lem:reduction-on-sphere}, we have
\begin{eqnarray}
\label{eqn:s-chi-lambda-diag-part-reduced-P-2-solo-fourier-heisenberg}
\lefteqn{S_{\chi \,e^{-i\lambda \,(\cdot)}} (x,x)}\\
&\sim&\frac{\lambda^{2\mathrm{d}+2}}{(2\pi)^{2\mathrm{d}+1}}\,\int _{-a}^a
\int_{B_{2\mathrm{d}}(\mathbf{0},\delta)}\int _0^{+\infty}\int_{-\epsilon}^\epsilon \int _{1/D}^D
\int_{S_+}
e^{i\,\lambda\,\Phi_2}\,r^{2\mathrm{d}}\,G\,d\theta\,d\mathbf{v}\,dt\,d\tau\,dr\,d\omega,\nonumber
\end{eqnarray}
where $\Phi_2$ is as in (\ref{eqn:phase-phi-heis}), and
\begin{equation}
G=:A\Big(x,x+(\theta,\mathbf{v}),\lambda t,\tau,\lambda r\omega\Big)\, F(r)\,\beta_+(\omega)\,
\mathcal{V}(\theta,\mathbf{v}).
\end{equation}

Our next task is to show that integration in $dt$ in (\ref{eqn:s-chi-lambda-diag-part-reduced-P-2-solo-fourier-heisenberg})
may also be restricted to a suitable compact interval in $(0,+\infty)$.

To this end, suppose $D'\gg 0$ and let $\rho :\mathbb{R}\rightarrow
\mathbb{R}$ be a $\mathcal{C}^\infty$ function satisfying $\rho \ge
0$, $\rho (t)=1$ if $|t|<2\, D'$, $\rho (t)=0$ if $|t|>3\,D'$.
Using the identity $G=\rho\,G+(1-\rho)\,G$ in
(\ref{eqn:s-chi-lambda-diag-part-reduced-P-2-solo-fourier-heisenberg}),
we get
\begin{equation}
\label{eqn:s-chi-lambda-diag-part-reduced-P-2-riscalato-tagliato-t}
S_{\chi \,e^{-i\lambda \,(\cdot)}} (x,x)\sim S_{\chi \,e^{-i\lambda
\,(\cdot)}}^{(1)} (x,x)+S_{\chi \,e^{-i\lambda \,(\cdot)}}^{(2)}
(x,x),
\end{equation}
where in $S_{\chi \,e^{-i\lambda \,(\cdot)}}^{(1)} (x,x)$
(respectively, $S_{\chi \,e^{-i\lambda \,(\cdot)}}^{(2)} (x,x)$) $G$ is replaced by
$\rho \,G$ (respectively, by $(1-\rho)\,G$).
Thus, $\int _0^{+\infty}dt$ may be replaced by $\int _0^{3D'} dt$ in
$S_{\chi \,e^{-i\lambda \,(\cdot)}}^{(1)} (x,x)$, and by $\int
_{2D'}^{+\infty}dt$ in $S_{\chi \,e^{-i\lambda \,(\cdot)}}^{(2)} (x,x)$,
respectively.

Recall that $B_{2\mathrm{d}}(\mathbf{0},\delta)$ is the image
of the preferred local chart on $M$ centered at $\pi(x)$
that underlies the chosen Heisenberg local chart on $X$.
Thus $\|\mathbf{v}\|<\delta$ for any $z=x+(\theta,\mathbf{v})\in X_1$.

\begin{lem}
\label{lem:second-term-S-decays}
If $0<\delta\ll 1$ and $D'\gg 0$,
then
$
S_{\chi\,e^{-i\lambda \,(\cdot)}}^{(2)} (x,x)
=O\left(\lambda^{-\infty}\right)
$
as $\lambda\rightarrow +\infty$.
\end{lem}

\textit{Proof.} By construction,
$d_{(x,x)}\psi=\left(\alpha_x,-\alpha_x\right)$ for any $x\in X$
\cite{bs}, \cite{z}, and more generally, $
d_{(e^{i\theta}x,x)}\psi=\left(e^{i\theta}\,\alpha_{e^{i\theta}x},-e^{-i\theta}\,\alpha_{x}\right),
$
for any $x\in X$ and
$e^{i\theta}\in S^1$.

By construction, $(e^{i\theta}x,x)=\big(x+(\theta,\mathbf{0}),x\big)$ in local coordinates.
Thus, if
$\Upsilon:(\theta,\mathbf{v})\mapsto
\psi\big(x+\left(\theta,\mathbf{v}\right),x\big)$,
then $\partial
_\theta\Upsilon\big(\theta,\mathbf{0}\big)=
\,e^{i\theta}$ for every $\theta$. Hence, if $\delta$ is
sufficiently small then
\begin{equation}
\label{eqn:bound-on-dpsi}
2\ge \left|\partial
_\theta\Upsilon\big(\theta,\mathbf{v}\big)\right
| \ge \frac 12,
\end{equation}
for every $z=x+(\theta,\mathbf{v})\in X_1$.

Suppose $D'>4D$; thus, $t\ge 2D'>8D$ on $\mathrm{supp}\{1-\rho\}$.
In view of (\ref{eqn:phase-phi-heis}) and the lower bound in
(\ref{eqn:bound-on-dpsi}), since $|\omega_0|\le 1$ and $r\le D$ we have
\begin{eqnarray}
\Big|\partial _\theta
\Phi_2\big(x,\theta,\mathbf{v},t,\tau,r,\omega\big)\Big|&=&\left|t\partial
_\theta\Upsilon\Big(\theta,\mathbf{v}\Big)-r\,\omega_0\right|\\
&\ge &\frac{t}{2}-r\ge\frac{t}{4}+\left(\frac{t}{4}-D\right)\ge \frac{t}{4}+D\nonumber
\end{eqnarray}
for all $t\in \mathrm{supp}\{1-\rho\}$.
The statement follows by integrating by parts in $\theta$.

\hfill Q.E.D.

\medskip

Thus $\int _0^{+\infty}dt$ may be replaced by
$\int _0^{3D'}dt$ in (\ref{eqn:s-chi-lambda-diag-part-reduced-P-2-solo-fourier-heisenberg}),
and the amplitude $G$ by $\rho\,G$.
We can similarly discard the contribution to $\int _0^{+\infty}dt$
coming from $(0,\varepsilon)$, if $0<\varepsilon\ll 1$.
Let $\upsilon :\mathbb{R}\rightarrow \mathbb{R}$
be $\mathcal{C}^\infty$, non-negative, and such that
$\upsilon (t)=0$ for $t\le \varepsilon$, $\upsilon (t)=1$ for $t>2\varepsilon$.
Thus
$S_{\chi \,e^{-i\lambda \,(\cdot)}} (x,x)\sim S_{\chi \,e^{-i\lambda \,(\cdot)}} (x,x)'+S_{\chi \,e^{-i\lambda \,(\cdot)}} (x,x)''$,
where $\rho\,G$ has been replaced by $\upsilon \,\rho\,G$ in $S_{\chi \,e^{-i\lambda \,(\cdot)}} (x,x)'$, and
by $(1-\upsilon)\,\rho\,G$ in $S_{\chi \,e^{-i\lambda \,(\cdot)}} (x,x)''$, respectively.

\begin{lem}
\label{lem:reduction-in-t-away-from-0}
If $\varepsilon$ is sufficiently small, then
$S_{\chi \,e^{-i\lambda \,(\cdot)}} (x,x)''=O\left(\lambda^{-\infty}\right)$ for
$\lambda\rightarrow +\infty$.
\end{lem}

In the following, recall that $\mathrm{supp}(\chi)\subseteq (-\epsilon,\epsilon)$, and that $\epsilon$ has
been chosen suitably small.

\medskip

\textit{Proof.}
Given (\ref{eqn:phase-phi-heis}),
$\partial _\tau\Phi_2=r\cdot \Big [q\big(x,\omega\big)+O(\tau)\Big]-1$.
Let $M$ be the maximum reached by $q$ on the unit sphere bundle of
$X$, and suppose $r<1/(2M)$.
Then $\Big|\partial _\tau\Phi_2\Big|\ge 1-\frac{M+O(\epsilon)}{2M}\ge \frac 13$,
since $\epsilon$ is small.
Integration by parts in $\tau$ shows that only a rapidly decreasing
contribution in $\lambda$ is lost if we multiply the amplitude by a  cut-off in
$r$ so as to assume $r>1/(2M)$.

Now suppose that $\varepsilon<\delta_1/(16\,M)$, and argue as in Lemma \ref{lem:second-term-S-decays},
using the upper bound in
(\ref{eqn:bound-on-dpsi}), and the fact that $\omega_0>\delta_1$ on $S_+$: if $t<2\varepsilon$,
\begin{eqnarray*}
\Big|\partial _\theta
\Phi_2\big(x,\theta,\mathbf{v},t,\tau,r,\omega\big)\Big|&=&\left|t\partial
_\theta\Upsilon\Big(\theta,\mathbf{v}\Big)-r\,\omega_0\right|\\
&\ge &r\omega_0-2t\ge \frac{\delta_1}{2M}-4\varepsilon\ge \frac{\delta_1}{4M}.\nonumber
\end{eqnarray*}
The statement follows by integration by parts in $\theta$.

\hfill Q.E.D.

\medskip

On the upshot, for $\lambda\rightarrow +\infty$
we only lose a rapidly decreasing term in (\ref{eqn:s-chi-lambda-diag-part-reduced-P-2-solo-fourier-heisenberg})
by replacing $\int_0^{+\infty}dt$ by $\int_{1/C}^{C}dt$ for some $C\gg 0$, and
$G$ by $\zeta (t)\,G $, where $\zeta \in \mathcal{C}^\infty (\mathbb{R})$
is supported in $(1/C,C)$, and identically one on
$(2/C,C/2)$.

As a further reduction, let $\gamma :\mathbb{C}^\mathrm{d}\rightarrow \mathbb{R}$ be
$\mathcal{C}^\infty$, nonnegative,
such that $\gamma (\mathbf{v})=1$ if $\|\mathbf{v}\|<2$ and
$\gamma (\mathbf{v})=0$ if $\|v\|>3$. Let us set
$\gamma_\lambda(\mathbf{v})=:\gamma \left(\lambda^{1/3}\mathbf{v}\right)$
($\lambda >0,\,\mathbf{v}\in \mathbb{C}^\mathrm{d}$).
We obtain $S_{\chi \,e^{-i\lambda \,(\cdot)}} (x,x)\sim S_{\chi \,e^{-i\lambda \,(\cdot)}} (x,x)^{\prime}_1+
S_{\chi \,e^{-i\lambda \,(\cdot)}} (x,x)^{\prime}_2$, where the amplitude is $\gamma_\lambda \,G$
in $S_{\chi \,e^{-i\lambda \,(\cdot)}} (x,x)^{\prime}_1$,
and $(1-\gamma_\lambda )\,G$ in $S_{\chi \,e^{-i\lambda \,(\cdot)}} (x,x)^{\prime}_2$,
respectively.

\begin{lem}
\label{lem:riscalato-in-v-due-termini-decade}
$S_{\chi \,e^{-i\lambda \,(\cdot)}} (x,x)^{\prime}_2=O\left(\lambda^{-\infty}\right)$ for
$\lambda\rightarrow +\infty$.
\end{lem}

\textit{Proof.} Let $\mathrm{dist}_X$ and $\mathrm{dist}_M$ be the Riemannian distance functions on
$X$ and $M$, respectively. If $m=:\pi(x)$, let
$m+\mathbf{v}$ be the point in $M$ with preferred coordinates $\mathbf{v}$.
By construction \cite{sz}, $m+\mathbf{v}=\pi\big(x+(\theta,\mathbf{v})\big)$
for every $\theta,\,\mathbf{v}$.

If $\mathbf{v}\in \mathrm{supp}\{1-\gamma_\lambda\}\cap B_{2\mathrm{d}}(\mathbf{0},\delta)$, then
$\delta>\|\mathbf{v}\|>2\,\lambda^{-1/3}$. By definition, in the preferred local coordinates
centered at $m$ the tangent space $T_mM$ is unitarily identified with $\mathbb{C}^\mathrm{d}$.
Therefore,
if $\delta$ is sufficiently small we may assume that
$\mathrm{dist}_M(m+\mathbf{v},m)\ge \frac 12\,\|\mathbf{v}\|$,
$\forall\,\mathbf{v}\in B_{2\mathrm{d}}(\mathbf{0},\delta)$. On the upshot,
\begin{eqnarray}
\label{eqn:bound-on-distance}
\mathrm{dist}_X\big(x+(\theta,\mathbf{v}),x\big)\ge\mathrm{dist}_M\big(m+\mathbf{v},m\big)\ge
\lambda^{-1/3},
\end{eqnarray}
$\forall \,\mathbf{v}\in \mathrm{supp}\{1-\gamma_\lambda\}\cap B_{2\mathrm{d}}(\mathbf{0},\delta),\,
\theta\in (-a,a)$.

By (1.4) of \cite{bs},
$\Im \psi (x,y)\ge C\,\mathrm{dist}_X(x,y)^2$ $\forall\,x,y\in X$ and a suitable
constant $C>0$.
In view of (\ref{eqn:bound-on-distance}),
we conclude that on the same support
\begin{eqnarray}
\label{eqn:bound-on-d-t-psi}
\Big|\partial _t\Phi_2\big(x,\theta,\mathbf{v},t,\tau,r,\omega\big)\Big|&=&\Big |\psi \big(x+(\theta,\mathbf{v}),x\big)\Big|
\nonumber \\
&\ge& \Im \psi \big(x+(\theta,\mathbf{v}),x\big)\ge C\lambda ^{-2/3}.
\end{eqnarray}
The statement then follows by integrating by parts in $t$, since each integration introduces a factor $\frac 1\lambda\cdot
\lambda ^{2/3}=\lambda ^{-1/3}$; notice that $t$ appears in the amplitude by a factor which has an asymptotic expansion
in $\lambda t$, with leading terms $\lambda ^\mathrm{d}t^\mathrm{d}$. Therefore, every
$t$-derivative of the amplitude
remains $O\left(\lambda ^\mathrm{d}\right)$.

\hfill Q.E.D.

\medskip

Lemma \ref{lem:riscalato-in-v-due-termini-decade} means that we only lose a rapidly decreasing term in (\ref{eqn:s-chi-lambda-diag-part-reduced-P-2-solo-fourier-heisenberg}) if, in addition to the previous reductions,
$\int _{B_{2\mathrm{d}}(\mathbf{0},\delta)}d\mathbf{v}$ is replaced by
$\int _{B_{2\mathrm{d}}(\mathbf{0},3\lambda^{-1/3})}d\mathbf{v}$, and the amplitude is multiplied by
$\gamma_\lambda(\mathbf{v})$.

Let us introduce the change of variables $\mathbf{v}\rightsquigarrow \mathbf{v}/(r\sqrt{\lambda})$.
Since $D^{-1}\le r\le D$, integration in $d\mathbf{v}$ is now over a ball of radius $O\left(\lambda^{1/6}\right)$.
We obtain:
\begin{eqnarray}
\label{eqn:s-chi-lambda-diag-part-reduced-P-2-lots-of-reductions}
\lefteqn{S_{\chi \,e^{-i\lambda \,(\cdot)}} (x,x)}\\
&\sim&\frac{\lambda^{\mathrm{d}+2}}{(2\pi)^{2\mathrm{d}+1}}\,\int _{-a}^a
\int_{\mathbb{C}^\mathrm{d}}\int _{1/C}^{C}\int_{-\epsilon}^\epsilon \int _{1/D}^D
\int_{S_+}
e^{i\,\lambda\,\Phi_3}\,H\,d\theta\,d\mathbf{v}\,dt\,d\tau\,dr\,d\omega,\nonumber
\end{eqnarray}
where
\begin{eqnarray*}
\lefteqn{\Phi_3=:\Phi_2\left(x,\theta,\frac{\mathbf{v}}{r\sqrt{\lambda}},t,\tau,r,\omega\right)}\\
&=&-\frac{\mathbf{v}}{\sqrt{\lambda}}\cdot\omega_1
+t\psi\left(x+\left(\theta,\frac{\mathbf{v}}{r\sqrt{\lambda}}\right),x\right)-r\,\theta \,\omega_0
+r\tau\,q\big(x,\omega\big)-\tau+O\left(\tau^2\right)\,r
\end{eqnarray*}
and
$$
H=:\gamma\left(\frac{\mathbf{v}}{r\lambda^{1/6}}\right)\,
\zeta(t)\,G\left(x,\theta,\frac{\mathbf{v}}{r\sqrt{\lambda}},t,\tau,r,\omega\right).
$$
In view of (65) in \cite{sz}, we have
\begin{eqnarray}
\lefteqn{t\psi\left(x+\left(\theta,\frac{\mathbf{v}}{r\sqrt{\lambda}}\right),x\right)}\nonumber\\
&=&
it\,\left[1-e^{i\theta}\right]+\left[\frac{it}{2r^2\lambda}\,\|\mathbf{v}\|^2\,e^{i\theta}
+t\,R_3^\psi\left(\frac{\mathbf{v}}{r\sqrt{\lambda}}\right)\,e^{i\theta}\right],
\end{eqnarray}
where $R_3^\psi$ vanishes to third order at the origin.

The latter term gives rise to a bounded exponential for $\|\mathbf{v}\|\lesssim \lambda ^{1/6}$;
as in \cite{sz}, this exponential may then be incorporated as part of the amplitude in an appopriate oscillatory
integral.
More precisely, we have:
\begin{eqnarray}
\label{eqn:s-chi-lambda-diag-part-reduced-P-2-lots-and-lots-of-reductions}
\lefteqn{S_{\chi \,e^{-i\lambda \,(\cdot)}} (x,x)}\\
&\sim&\frac{\lambda^{\mathrm{d}+2}}{(2\pi)^{2\mathrm{d}+1}}\,\int_{\mathbb{C}^\mathrm{d}}\int_{S_+}
e^{-i\sqrt{\lambda}\,\mathbf{v}\cdot \omega_1}\cdot \left[
\int _{-a}^a
\int _{1/C}^{C}\int_{-\epsilon}^\epsilon \int _{1/D}^D
e^{i\,\lambda\,\Psi}\,\widetilde{H}\,d\theta\,dt\,d\tau\,dr\right]\,d\mathbf{v}\,d\omega,\nonumber
\end{eqnarray}
where
\begin{equation}
\Psi=:it\,\left[1-e^{i\theta}\right]-r\,\theta \,\omega_0
+r\tau\,q\big(x,\omega\big)-\tau+O\left(\tau^2\right)\,r,
\end{equation}
\begin{equation}
\widetilde{H}=:e^{-(t/2r^2)\,\|\mathbf{v}\|^2\,e^{i\theta}
+i\lambda\,t\,R_3^\psi\left(\frac{\mathbf{v}}{r\sqrt{\lambda}}\right)\,e^{i\theta}}\,H.
\end{equation}

Let us first consider the asymptotics of the inner integral in
(\ref{eqn:s-chi-lambda-diag-part-reduced-P-2-lots-and-lots-of-reductions}),
an oscillatory integral $d\theta\,dt\,d\tau\,dr$ with phase $\Psi$, and amplitude
$\widetilde{H}$ that may be expanded in descending powers of $\lambda^{-1/2}$
for $\lambda\rightarrow +\infty$; here
$\mathbf{v}\in \mathbb{C}^\mathrm{d}$, $\omega\in S_+$ play the role of
parameters.
Using the the ellipticity of $q$ and the fact that $|\tau|<\epsilon$, $\epsilon$ small,
a straightforward computation yields the following.

\begin{lem}
\label{lem:stationary-point}
On the given domain of integration, $\Psi$ has a unique stationary point
given by
$(\theta_0,t_0,\tau_0,r_0)=\big(0,\omega_0/q(x,\omega),0,1/q(x,\omega)\big)$.
Furthermore, if $\Psi''$ is the Hessian of $\Psi$ at this critical point then
$$
\det\left(\frac{\lambda \Psi''}{2\pi\,i}\right)=\left(\frac{\lambda}{2\pi}\right)^4\,q(x,\omega)^2.
$$
\end{lem}

We are thus in a position to apply the stationary phase Lemma, following \S 5 of \cite{sz},
and conclude that for any integer $N\gg 0$ the inner integral in
(\ref{eqn:s-chi-lambda-diag-part-reduced-P-2-lots-and-lots-of-reductions}) is given by
$S_N(x,\omega,\mathbf{v})+R_N(x,\omega,\mathbf{v})$, where
\begin{eqnarray}
\label{eqn:inner-integral-expanded}
\lefteqn{S_N(x,\omega,\mathbf{v})=\lambda^{\mathrm{d}-2}\,(2\pi)^{2}\,\chi(0)}\\
&&\cdot \frac{\omega_0^\mathrm{d}}{q(x,\omega)^{\mathrm{d}+1}}\,\,e^{-\frac 12\omega_0\,
q(x,\omega)\,\|\mathbf{v}\|^2}\beta_+(\omega)\,s(x,x)\,\left(1+\sum _{j=1}^NF_j(\omega,\mathbf{v})
\lambda^{-j/2}\right),\nonumber
\end{eqnarray}
where $F_j$ is polynomial in $\mathbf{v}$ and
$\big|R_N(x,\omega,\mathbf{v})\big|\le C_N\,\lambda^{\mathrm{d}-2-(N+1)/2}\,e^{-a\,\|\mathbf{v}\|^2}$, for some $a>0$.
Therefore, the asymptotic expansion may be integrated term by term in $d\mathbf{v}\,d\omega$.

In order to apply the stationary phase Lemma again, we set $\mu=\sqrt{\lambda}$. Since
$s_0(x,x)=\pi^{-\mathrm{d}}$, we have
\begin{eqnarray}
\label{eqn:s-chi-lambda-diag-part-reduced-P-2-stationary-phase-mu}
\lefteqn{S_{\chi \,e^{-i\lambda \,(\cdot)}} (x,x)}\\
&\sim&\frac{2\pi}{\pi^\mathrm{d}} \,\left(\frac{\mu^2}{2\pi}\right)^{2\mathrm{d}}\,\chi (0)\cdot \left(
\int_{\mathbb{C}^\mathrm{d}}\int_{S_+} e^{-i\mu\,\mathbf{v}\cdot \omega_1-\frac 12\omega_0\,
q(x,\omega)\,\|\mathbf{v}\|^2} \,\frac{\omega_0^\mathrm{d}}{q(x,\omega)^{\mathrm{d}+1}}
\,\beta_+(\omega)\,d\mathbf{v}\,d\omega \right.   \nonumber\\
&&\left.+\sum _{j=1}^{+\infty}\mu^{-j}\int_{\mathbb{C}^\mathrm{d}}\int_{S_+}
e^{-i\mu\,\mathbf{v}\cdot \omega_1-\frac 12\omega_0\,
q(x,\omega)\,\|\mathbf{v}\|^2} \,\frac{\omega_0^\mathrm{d}}{q(x,\omega)^{\mathrm{d}+1}}
\,\beta_+(\omega)\,F_j(\omega,\mathbf{v})
\,d\mathbf{v}\,d\omega \right).\nonumber
\end{eqnarray}

To obtain an asymptotic expansion for $S_{\chi \,e^{-i\lambda \,(\cdot)}} (x,x)$ as $\lambda\rightarrow +\infty$,
we are thus reduced to
determining asymptotic expansions for each of
the summands in (\ref{eqn:s-chi-lambda-diag-part-reduced-P-2-stationary-phase-mu}) as $\mu\rightarrow +\infty$.
To this end, we remark that $\omega_1$ is a system of local coordinates on $S_+\subseteq S^{2\mathrm{d}}$, and
$\omega_0=\left(1-\|\omega_1\|^2\right)^{1/2}$.
Each summand in (\ref{eqn:s-chi-lambda-diag-part-reduced-P-2-stationary-phase-mu}) is an oscilatory integral whose phase is the quadratic form $(\mathbf{v},\omega_1)\mapsto -\mathbf{v}\cdot \omega_1$
on $\mathbb{R}^{\mathrm{4d}}\cong \mathbb{R}^{2\mathrm{d}}\times \mathbb{R}^{2\mathrm{d}}$.
The origin $\mathbf{v}=\mathbf{0},\,\omega_1=\mathbf{0}$ is the only stationary point, and the Hessian has determinant
one. Thus, the stationary phase Lemma applies and as $\mu\rightarrow +\infty$ we are left with an asymptotic expansion
\begin{eqnarray}
\label{eqn:asympt-exp-1st-summand}
\lefteqn{S_{\chi \,e^{-i\lambda \,(\cdot)}} (x,x)}\\
&\sim&
\frac{2\pi}{\pi^\mathrm{d}} \,\left(\frac{\mu^2}{2\pi}\right)^{2\mathrm{d}}\,\chi (0)\cdot \left(\frac{2\pi}{\mu}\right)
^{2\mathrm{d}}\,q\big(x,(1,\mathbf{0})\big)^{-(\mathrm{d}+1)}\cdot \left(1+\sum_{j\ge 1}c_j\,\mu^{-j}\right)\nonumber\\
&=&2\pi\,\left(\frac{\lambda}{\pi}\right)^{\mathrm{d}}\,\chi (0)\,
\varsigma_T(x)^{-(\mathrm{d}+1)}\cdot \left(1+\sum_{j\ge 1}c_j\,\lambda^{-j/2}\right).\nonumber
\end{eqnarray}

\begin{lem}
\label{lem:c1vanishes}
$c_1=0$.
\end{lem}

\textit{Proof.} $c_1$ is the coefficient of $\mu^{2\mathrm{d}-1}$ in the second line of (\ref{eqn:asympt-exp-1st-summand}).
Looking at (\ref{eqn:s-chi-lambda-diag-part-reduced-P-2-stationary-phase-mu}), there are two possible contributions to
this coefficient.

One is the second term in the asymptotic expansion of the first summand in (\ref{eqn:s-chi-lambda-diag-part-reduced-P-2-stationary-phase-mu}); by the stationary phase Lemma, this is multiple of
$$
\langle D_\mathbf{v},D_{\omega_1}\rangle \left(e^{-\frac 12\omega_0\,
q(x,\omega)\,\|\mathbf{v}\|^2} \,\frac{\omega_0^\mathrm{d}}{q(x,\omega)^{\mathrm{d}+1}}\,\beta_+(\omega)\right)(\mathbf{0},\mathbf{0})=0.
$$

The other possible contribution comes from the leading term in the asymptotic expansion of the second
summand in (\ref{eqn:s-chi-lambda-diag-part-reduced-P-2-stationary-phase-mu}), the one with $j=1$;
hence it is a multiple of the
evaluation at the critical point of the amplitude of that oscillatory integral. On the other hand, the appearance
of the fractional powers of $\lambda^{-1}$ in the amplitude originates from the Taylor expansion of $\widetilde{H}$
in $\mathbf{v}/(r\sqrt{\lambda})$; hence $F_1$ is linear in $\mathbf{v}$, and so it vanishes at the critical point.

\hfill Q.E.D.

\begin{cor}
\label{cor:summingup}
There exists $\epsilon>0$ such that for every $\chi \in \mathcal{C}^\infty_0\big((-\epsilon,\epsilon)\big)$

$$
S_{\chi \,e^{-i\lambda \,(\cdot)}} (x,x)=\left\{
\begin{array}{lcc}
  2\pi\,\left(\frac{\lambda}{\pi}\right)^{\mathrm{d}}\,\chi (0)\,\varsigma_T(x)^{-(\mathrm{d}+1)}
+O\left(\lambda^{\mathrm{d}-1}\right) & \mathrm{for} & \lambda\rightarrow +\infty \\
   &  &  \\
  O\left(\lambda^{-\infty}\right) & \mathrm{for} & \lambda\rightarrow -\infty
\end{array}\right.
.
$$
\end{cor}

To complete the proof of Theorem \ref{thm:local-weyl-law}, we shall now closely follow \S 12 of \cite{gs},
hence the exposition will be rather sketchy.

Le us fix once and forall $\epsilon>0$ such that the conclusion of Corollary \ref{cor:summingup}
holds true. One can see that $S_{\chi \,e^{-i\lambda \,(\cdot)}}$ is a smoothing operator, with kernel
$$
S_{\chi \,e^{-i\lambda \,(\cdot)}} (x,y)=\sum _j\widehat{\chi}(\lambda-\lambda_j)\,e_j(x)\,\overline{e_j(y)}
\,\,\,\,\,\,\,\,(x,y\in X),
$$
hence the Corollary implies that for $\lambda\rightarrow +\infty$
$$
\sum _j\widehat{\chi}(\lambda-\lambda_j)\,\big|e_j(x)\big|^2=
2\pi\,\left(\frac{\lambda}{\pi}\right)^{\mathrm{d}}\,\chi (0)\,\varsigma_T(x)^{-(\mathrm{d}+1)}
+O\left(\lambda^{\mathrm{d}-1}\right).
$$

Applying this to  $\chi \in \mathcal{C}^\infty_0\big((-\epsilon,\epsilon)\big)$ such that
$\chi (0)=1$ and $\widehat{\chi}\ge 0$ shows the following:

\begin{cor}
\label{cor:incremental-estimate}
There exists a constant
$E>0$ such that
$\mathcal{T}(\lambda+1,x)-\mathcal{T}(\lambda,x)\le E\,\lambda^{\mathrm{d}}$ for $\lambda\rightarrow +\infty$,
whence
$\mathcal{T}(\lambda,x)=O\left(\lambda^{\mathrm{d}+1}\right)$. These estimates are uniform in $x\in X$.
\end{cor}

Let us define measures $d\mathcal{T}_x=:\sum _j\big|e_j(x)\big|^2\,\delta_{\lambda_j}$ ($x\in X$)
on $\mathbb{R}$, so that
\begin{equation}
\label{eqn:measure-integrated}
\mathcal{T}(\lambda,x)=\int_{-\infty}^\lambda d\mathcal{T}_x(\eta),\,\,\,\,\,\,\,\,
S_{\chi \,e^{-i\lambda \,(\cdot)}} (x,x)=\int_{-\infty}^{+\infty}\widehat{\chi}(\lambda-\eta)\,d\mathcal{T}_x(\eta).
\end{equation}

Under the same assumption on $\chi$, set $G(\lambda)=:\int_{-\infty}^\lambda \widehat{\chi}(\tau)\,d\tau$; the integral
$\int _{-\infty}^{+\infty}G(\lambda-\eta)\,d\mathcal{T}_x(\eta)$ may be computed in two differen manners.

On the hand, we use the change of variable $\tau \rightsquigarrow \tau-\eta$ and the Tonelli-Fubini theorem
(recall that $\widehat{\chi}\ge 0$) to change the order of integration and obtain
$\int _{-\infty}^{+\infty}G(\lambda-\eta)\,d\mathcal{T}_x(\eta)=
\int_{-\infty}^\lambda S_{\chi \,e^{-i\tau \,(\cdot)}} (x,x)\,d\tau$. Using Corollary \ref{cor:summingup},
we get
\begin{equation}
\label{eqn:1stcomputation-of-int}
\int _{-\infty}^{+\infty}G(\lambda-\eta)\,d\mathcal{T}_x(\eta)
=\frac{2\,\pi^2}{(\mathrm{d}+1)}\cdot
\left(\frac{\lambda}{\pi\,\varsigma_T(x)}\right)^{\mathrm{d}+1}
+O\left(\lambda^{d}\right).
\end{equation}

On the other hand, if $H$ is the Heaviside function, we have the chain of equalities:
\begin{eqnarray}
\label{eqn:2ndcomputation-of-int}
\lefteqn{\int _{-\infty}^{+\infty}G(\lambda-\eta)\,d\mathcal{T}_x(\eta)=\sum_jG(\lambda-\lambda_j)\,\big|e_j(x)\big|^2}\\
&=&\sum _j\left(\int _{-\infty}^{\lambda-\lambda_j}\widehat{\chi}(\tau)\,d\tau\right)\cdot \big|e_j(x)\big|^2
=\sum_j\int _{-\infty}^{+\infty}H(\lambda-\lambda_j-\tau)\,\widehat{\chi}(\tau)
\cdot \big|e_j(x)\big|^2\,d\tau\nonumber\\
&=&\int _{-\infty}^{+\infty}\left(\sum_jH(\lambda-\lambda_j-\tau)\,
\, \big|e_j(x)\big|^2\right)\cdot \widehat{\chi}(\tau)\,d\tau =
\int _{-\infty}^{+\infty}\mathcal{T}(\lambda-\tau,x)\cdot \widehat{\chi}(\tau)\,d\tau \nonumber\\
&=&\mathcal{T}(\lambda,x)\cdot \int _{-\infty}^{+\infty}\widehat{\chi}(\tau)\,d\tau +
\int _{-\infty}^{+\infty}\Big(\mathcal{T}(\lambda-\tau,x)-\mathcal{T}(\lambda,x)\Big)\cdot \widehat{\chi}(\tau)\,d\tau\nonumber\\
&=&2\pi\,\mathcal{T}(\lambda,x)+O\left(\lambda^{\mathrm{d}}\right),\nonumber
\end{eqnarray}
where the estimate on the second summand is a consequence of Corollary \ref{cor:incremental-estimate} and the fact that
$\widehat{\chi}$ is of rapid decay.

Theorem \ref{thm:local-weyl-law} follows from (\ref{eqn:1stcomputation-of-int}) and
(\ref{eqn:2ndcomputation-of-int}).

\hfill Q.E.D.

\section{Proof of Corollary \ref{cor:global-weyl-law}}

In order to prove Corollary \ref{cor:global-weyl-law}, let us dwell on the symplectic structure
$\Omega_\Sigma$ of the cone $\Sigma$
induced by the canonical symplectic structure of $T^*X$.

The action $r$ of $S^1$ on $X$ lifts to an action on $\Sigma$, given by
$e^{i\theta}\cdot (x,\alpha_x)=:\left(e^{i\theta}\cdot x,\alpha_{e^{i\theta}\cdot x}\right)$;
denote by $\frac{\partial}{\partial \theta}$ the generator of this action.

$\mathbb{R}_+=(0,+\infty)$ also acts on the cone $\Sigma$, by $r\cdot (x,s\,\alpha_x)=:(x,rs\,\alpha_x)$;
let $\frac{\partial}{\partial r}$ be the generator of the latter action; it spans the vertical tangent bundle
of the projection  $p:\Sigma\rightarrow X$.

Thus $\frac{\partial}{\partial \theta}$ and $\frac{\partial}{\partial r}$
are smooth nowhere vanishing vector fields on $\Sigma$, and they are actually everywhere linearly independent.
Let $\mathfrak{V}\subseteq T\Sigma$ be the rank-2 vector subbundle of the tangent bundle of $\Sigma$
spanned by $\frac{\partial}{\partial \theta}$ and $\frac{\partial}{\partial r}$. On $\mathfrak{V}$ let us consider
the symplectic structure $\Omega _\mathfrak{V}$ uniquely determined by
$\Omega _\mathfrak{V}\left( \frac{\partial}{\partial \theta},\frac{\partial}{\partial r}\right)=1$.
%In Heisenberg local coordinates, say, we have $\Omega _\mathfrak{V}=d\theta\wedge dr$.

For every $r>0$, the map
$\sigma_r:X\rightarrow \Sigma$, $x\mapsto (x,r\alpha_x)$, is a smooth section of $p$, smoothly varying with $r$.
Thus $d_x\sigma_r\big(T_xX\big)\subseteq T_{(x,r\alpha_x)}\Sigma$ is the fiber of a vector subbundle of rank
$2\mathrm{d}+1$ of $T\Sigma$;
therefore, $\mathfrak{H}_{(x,r\alpha_x)}=:d_x\sigma_r\big(\ker(\alpha_x)\big)$
is the fiber of a vector subbundle of rank $2\mathrm{d}$ of $T\Sigma$, isomorphic to the pull-back $q^*\big(TM\big)$,
where $q=:\pi\circ p:\Sigma\rightarrow M$. We have $T\Sigma=\mathfrak{H}\oplus \mathfrak{V}$.

On $\mathfrak{H}\cong q^*(TM)$ we have the symplectic structure $q^*(\omega)$; let $\Omega_\mathfrak{H}$ be the symplectic
structure on $\mathfrak{H}$ defined by
$\left(\Omega_\mathfrak{H}\right)_{(x,r\alpha_x)}=2\,r\,q^*(\omega)_{(x,r\alpha_x)}$.

\begin{lem}
\label{lem:symp-str-cone}
$\big(T\Sigma,\Omega_\Sigma)\cong \big(\mathfrak{H},\Omega_\mathfrak{H}\big)\oplus \big(\mathfrak{V},\Omega_\mathfrak{V}\big)$
as symplectic vector bundles.
\end{lem}

\textit{Proof.} Suppose $x\in X$, and let $U\subseteq X$ be an open neighborhood of $x$,
on which a local coordinate chart $\mathbf{t}=(t_i):U\rightarrow B$ is defined;
here $B$ is an open ball in $\mathbb{R}^{2\mathrm{d}+1}$, and $\mathbf{t}$ is a diffeomorfism.

Using the basis $\{dt_i\}_i$,
a differential 1-form $\widetilde{a}=\sum_i a_i\,dt_i$ on $U$ may be represented
as a smooth function $a:B\rightarrow \mathbb{R}^{2\mathrm{d}+1}$, $\mathbf{t}\mapsto
\big(a_i(\mathbf{t})\big)$.
%;if $g:U\rightarrow \mathbb{R}^{2\mathrm{d}+1}$ is smooth, let $\widetilde{g}=\sum g_i\,dt_i$ be the corresponding
%1-form.

Using the basis $\{dt_i\wedge dt_j\}_{i<j}$,
a differential 2-form $\widehat{\nu}=\sum _{i<j}\nu _{ij}\,dt_i\wedge dt_j$ may be represented
as a smooth function $\nu:B\rightarrow A_{2\mathrm{d}+1}$,
$y\mapsto \nu(\mathbf{t})=\big[\nu_{ij}(\mathbf{t})\big]$; here
$A_{2\mathrm{d}+1}$ is the vector space of skew-symmetric $(2\mathrm{d}+1)\times (2\mathrm{d}+1)$ real matrices.
For any $y\in U$, $\nu(y)$ is the matrix representing $\widehat{\nu}_y$
as a skew-symmetric bilinear pairing on $T_yX$.
%;if $S:U\rightarrow A_{2\mathrm{d}+1}$ is smooth, let $\widehat{S}=\sum _{i<j}S _{ij}\,dt_i\wedge dt_j$
%be the corresponding 2-form.

If $g:B\rightarrow \mathbb{R}^{2\mathrm{d}+1}$ is smooth, let $\mathrm{Jac}(g)$ be the its Jacobian matrix;
then $d\widetilde{g}=\widehat{S(g)}$, where $S(g)=:\mathrm{Jac}(g)^t-\mathrm{Jac}(g)$.

If $\alpha =\widetilde{a}$ on $U$, and $y\in X$ has local coordinates $\mathbf{t}\in B$,
the horizontal subspace $\ker (\alpha_y)\subseteq T_yX$
corresponds to the hyperplane $\big\{\mathbf{v}:a(\mathbf{t})\cdot \mathbf{v}
=0\big\}\subseteq \mathbb{R}^{2\mathrm{d}+1}$. Let $\kappa (\mathbf{t})\in \mathbb{R}^{2\mathrm{d}+1}$
correspond to $\left.\frac{\partial}{\partial \theta}\right|_y$.

On $\mathbb{R}^{2\mathrm{d}+1}\times \mathbb{R}^{2\mathrm{d}+1}$, with linear coordinates $(\mathbf{q},\mathbf{p})$,
consider the symplectic structure $d\mathbf{p}\wedge d\mathbf{q}$; the local chart $\mathbf{t}$ determines
a the symplectic coordinate chart $T^*(U)\cong  B\times \mathbb{R}^{2\mathrm{d}+1}$.
Locally in this chart, $\Sigma$ is parametrized as $\big\{\big(\mathbf{t},r\,a(\mathbf{t})\big):\mathbf{t}\in B,\,r>0\big\}$.
Therefore, if $r>0$ and
$z=:(y,r\,\alpha_y)$ then the tangent space $T_{z}\Sigma$
corresponds to the vector subspace
\begin{eqnarray}
\label{eqn:tangent-space-sigma}
\lefteqn{\left\{\left(
                  \begin{array}{c}
                    \mathbf{v} \\
                    r\,\mathrm{Jac}_\mathbf{t}(a)\,\mathbf{v} \\
                  \end{array}
                \right):\mathbf{v}\in \mathbb{R}^{2\mathrm{d}+1}\right\}\oplus \mathrm{span}\left\{\left(
                  \begin{array}{c}
                    0 \\
                    a(\mathbf{t}) \\
                  \end{array}
                \right)\right\}}\\
                &=&\left\{\left(
                  \begin{array}{c}
                    \mathbf{v} \\
                    r\,\mathrm{Jac}_\mathbf{t}(a)\,\mathbf{v} \\
                  \end{array}
                \right):\mathbf{v}\in \mathbb{R}^{2\mathrm{d}+1},\,a(\mathbf{v})=0\right\}\nonumber \\
                &&\oplus
                \mathrm{span}\left\{\left(
                  \begin{array}{c}
                    \kappa (\mathbf{t})\\
                    r\,\mathrm{Jac}_\mathbf{t}(a)\,\kappa (\mathbf{t})\\
                  \end{array}
                \right),\,\left(
                  \begin{array}{c}
                    0 \\
                    a(\mathbf{t}) \\
                  \end{array}
                \right)\right\}.\nonumber
                \end{eqnarray}

Since the chart is symplectic, the statement follows from (\ref{eqn:tangent-space-sigma}) recalling that
$d\alpha=2\,\omega$.

\hfill Q.E.D.

\medskip

We can now prove Corollary \ref{cor:global-weyl-law}.
Given Lemma \ref{lem:symp-str-cone}, the volume form on $\Sigma$ is
\begin{eqnarray}
\label{eqn:vol-form-cone}
\lefteqn{dV_\Sigma=\frac{1}{(d+1)!}\,\Omega_\Sigma^{\wedge (\mathrm{d}+1)}=
\frac{1}{\mathrm{d!}}\,\Omega_\mathfrak{H}^{\wedge \mathrm{d}}
\wedge \Omega_\mathfrak{V}}\\
&=&2^{\mathrm{d}}\,r^{\mathrm{d}}\,q^*\big(dV_M\big)\wedge \big(d\theta\wedge dr\big)
=2^{\mathrm{d}+1}\,\pi\,r^{\mathrm{d}}\,p^*(d\mu_X\big)\wedge  dr.\nonumber
\end{eqnarray}

By (\ref{eqn:vol-form-cone}), the symplectic volume of $\Sigma_1=\big\{(x,r\,\alpha_x)\in \Sigma: r\le 1/\varsigma_T(x)\big\}$
is
\begin{eqnarray}
\mathrm{Vol}(\Sigma_1)&=&2^{\mathrm{d}+1}\,\pi\,\int _X\left(\int _0^{1/\varsigma_T(x)}r^{\mathrm{d}}\,dr\right)\,d\mu_X\\
&=&2^{\mathrm{d}+1}\,\frac{\pi}{\mathrm{d}+1}\,\int _X\left(\frac{1}{\varsigma_T(x)}\right)^{\mathrm{d}+1}\,d\mu_X.\nonumber
\end{eqnarray}

Given this, Theorem \ref{thm:local-weyl-law} and (\ref{eqn:counting-trace}) imply
\begin{eqnarray}
N_T(\lambda)&=&\left(\frac{\lambda}{\pi}\right)^{\mathrm{d}+1}\,\frac{\pi}{d+1}\,
\int_X\left(\frac{1}{\varsigma_T(x)}\right)^{\mathrm{d}+1}\,d\mu_X
+O\left(\lambda^\mathrm{d}\right)\\
&=&\left(\frac{\lambda}{2\pi}\right)^{\mathrm{d}+1}\,\mathrm{Vol}(\Sigma_1)
+O\left(\lambda^\mathrm{d}\right).\nonumber
\end{eqnarray}

\hfill Q.E.D.

\end{document}